\documentclass[11pt]{article}
\usepackage{times,amsmath,amssymb,graphicx,theorem}
\usepackage{pstricks,pst-node,pst-text,pst-tree,pstricks-add}
\makeatletter
\def\section{\@startsection{section}{1}{0pt}{-3.25ex plus -1ex minus 
-.2ex}{1.5ex plus .2ex minus .3ex}{\normalfont\large\bf}}
\renewcommand\subsection{\@startsection{subsection}{2}{\z@}%
                                     {-3.25ex\@plus -1ex \@minus -.2ex}%
                                     {1.5ex \@plus .2ex}%
                                     {\normalfont\normalsize\bfseries}}
\makeatother
\renewenvironment{abstract}
{\vspace*{-1.8ex}\begin{quotation}\small\noindent\textbf{Abstract. }}{\end{quotation}}
\newcommand{\parag}[1]{\vspace*{-1ex}\paragraph{#1.}\!\!\!}
\newcommand{\defn}[1]{{\textit{\textbf{#1}}}}

\theoremheaderfont{\scshape}
\theorembodyfont{\normalfont\slshape}
\theoremstyle{plain}

\newtheorem{proposition}{Proposition}
\newtheorem{corollary}{Corollary}
\newtheorem{theorem}{Theorem}
\newenvironment{proof}{\begin{trivlist}\item{}\normalfont\textit{Proof.}}{\hfill$\square$\end{trivlist}}

\newcommand{\nth}{^{\text{th}}}

\newcommand{\ie}{\emph{i.e.}}
\newcommand{\eg}{\emph{e.g.}}
\newcommand{\cf}{\emph{cf.}}

\newdimen\arrayruleHwidth
\makeatletter
\def\Hline{\noalign{\ifnum0=`}\fi\hrule \@height \arrayruleHwidth
  \futurelet \@tempa\@xhline}
\makeatother

\newdimen\proofrulebreadth \proofrulebreadth=.05em
\newdimen\proofdotseparation \proofdotseparation=1.25ex
\newdimen\proofrulebaseline \proofrulebaseline=2ex
\newcount\proofdotnumber \proofdotnumber=3
\let\then\relax
\def\hfi{\hskip0pt plus.0001fil}
\mathchardef\squigto="3A3B
\newif\ifinsideprooftree\insideprooftreefalse
\newif\ifonleftofproofrule\onleftofproofrulefalse
\newif\ifproofdots\proofdotsfalse
\newif\ifdoubleproof\doubleprooffalse
\let\wereinproofbit\relax
\newdimen\shortenproofleft
\newdimen\shortenproofright
\newdimen\proofbelowshift
\newbox\proofabove
\newbox\proofbelow
\newbox\proofrulename
\def\shiftproofbelow{\let\next\relax\afterassignment\setshiftproofbelow\dimen0 }
\def\shiftproofbelowneg{\def\next{\multiply\dimen0 by-1 }%
\afterassignment\setshiftproofbelow\dimen0 }
\def\setshiftproofbelow{\next\proofbelowshift=\dimen0 }
\def\setproofrulebreadth{\proofrulebreadth}

\def\prooftree{
\ifnum  \lastpenalty=1
\then   \unpenalty
\else   \onleftofproofrulefalse
\fi
\ifonleftofproofrule
\else   \ifinsideprooftree
        \then   \hskip.5em plus1fil
        \fi
\fi
\bgroup
\setbox\proofbelow=\hbox{}\setbox\proofrulename=\hbox{}%
\let\justifies\proofover\let\leadsto\proofoverdots\let\Justifies\proofoverdbl
\let\using\proofusing\let\[\prooftree
\ifinsideprooftree\let\]\endprooftree\fi
\proofdotsfalse\doubleprooffalse
\let\thickness\setproofrulebreadth
\let\shiftright\shiftproofbelow \let\shift\shiftproofbelow
\let\shiftleft\shiftproofbelowneg
\let\ifwasinsideprooftree\ifinsideprooftree
\insideprooftreetrue
\setbox\proofabove=\hbox\bgroup$\displaystyle 
\let\wereinproofbit\prooftree
\shortenproofleft=0pt \shortenproofright=0pt \proofbelowshift=0pt
\onleftofproofruletrue\penalty1
}

\def\eproofbit{
\ifx    \wereinproofbit\prooftree
\then   \ifcase \lastpenalty
        \then   \shortenproofright=0pt  
        \or     \unpenalty\hfil         
        \or     \unpenalty\unskip       
        \else   \shortenproofright=0pt  
        \fi
\fi
\global\dimen0=\shortenproofleft
\global\dimen1=\shortenproofright
\global\dimen2=\proofrulebreadth
\global\dimen3=\proofbelowshift
\global\dimen4=\proofdotseparation
\global\count255=\proofdotnumber
$\egroup  
\shortenproofleft=\dimen0
\shortenproofright=\dimen1
\proofrulebreadth=\dimen2
\proofbelowshift=\dimen3
\proofdotseparation=\dimen4
\proofdotnumber=\count255
}

\def\proofover{
\eproofbit 
\setbox\proofbelow=\hbox\bgroup 
\let\wereinproofbit\proofover
$\displaystyle
}%
\def\proofoverdbl{
\eproofbit 
\doubleprooftrue
\setbox\proofbelow=\hbox\bgroup 
\let\wereinproofbit\proofoverdbl
$\displaystyle
}%
\def\proofoverdots{
\eproofbit 
\proofdotstrue
\setbox\proofbelow=\hbox\bgroup 
\let\wereinproofbit\proofoverdots
$\displaystyle
}%
\def\proofusing{
\eproofbit 
\setbox\proofrulename=\hbox\bgroup 
\let\wereinproofbit\proofusing
\kern0.3em$
}

\def\endprooftree{
\eproofbit 
  \dimen5 =0pt
\dimen0=\wd\proofabove \advance\dimen0-\shortenproofleft
\advance\dimen0-\shortenproofright
\dimen1=.5\dimen0 \advance\dimen1-.5\wd\proofbelow
\dimen4=\dimen1
\advance\dimen1\proofbelowshift \advance\dimen4-\proofbelowshift
\ifdim  \dimen1<0pt
\then   \advance\shortenproofleft\dimen1
        \advance\dimen0-\dimen1
        \dimen1=0pt
        \ifdim  \shortenproofleft<0pt
        \then   \setbox\proofabove=\hbox{%
                        \kern-\shortenproofleft\unhbox\proofabove}%
                \shortenproofleft=0pt
        \fi
\fi
\ifdim  \dimen4<0pt
\then   \advance\shortenproofright\dimen4
        \advance\dimen0-\dimen4
        \dimen4=0pt
\fi
\ifdim  \shortenproofright<\wd\proofrulename
\then   \shortenproofright=\wd\proofrulename
\fi
\dimen2=\shortenproofleft \advance\dimen2 by\dimen1
\dimen3=\shortenproofright\advance\dimen3 by\dimen4
\ifproofdots
\then
        \dimen6=\shortenproofleft \advance\dimen6 .5\dimen0
        \setbox1=\vbox to\proofdotseparation{\vss\hbox{$\cdot$}\vss}%
        \setbox0=\hbox{%
                \advance\dimen6-.5\wd1
                \kern\dimen6
                $\vcenter to\proofdotnumber\proofdotseparation
                        {\leaders\box1\vfill}$%
                \unhbox\proofrulename}%
\else   \dimen6=\fontdimen22\the\textfont2 
        \dimen7=\dimen6
        \advance\dimen6by.5\proofrulebreadth
        \advance\dimen7by-.5\proofrulebreadth
        \setbox0=\hbox{%
                \kern\shortenproofleft
                \ifdoubleproof
                \then   \hbox to\dimen0{%
                        $\mathsurround0pt\mathord=\mkern-6mu%
                        \cleaders\hbox{$\mkern-2mu=\mkern-2mu$}\hfill
                        \mkern-6mu\mathord=$}%
                \else   \vrule height\dimen6 depth-\dimen7 width\dimen0
                \fi
                \unhbox\proofrulename}%
        \ht0=\dimen6 \dp0=-\dimen7
\fi
\let\doll\relax
\ifwasinsideprooftree
\then   \let\VBOX\vbox
\else   \ifmmode\else$\let\doll=$\fi
        \let\VBOX\vcenter
\fi
\VBOX   {\baselineskip\proofrulebaseline \lineskip.2ex
        \expandafter\lineskiplimit\ifproofdots0ex\else-0.6ex\fi
        \hbox   spread\dimen5   {\hfi\unhbox\proofabove\hfi}%
        \hbox{\box0}%
        \hbox   {\kern\dimen2 \box\proofbelow}}\doll%
\global\dimen2=\dimen2
\global\dimen3=\dimen3
\egroup 
\ifonleftofproofrule
\then   \shortenproofleft=\dimen2
\fi
\shortenproofright=\dimen3
\onleftofproofrulefalse
\ifinsideprooftree
\then   \hskip.5em plus 1fil \penalty2
\fi
}

\newcount\PLv\newcount\PLw\newcount\PLx\newcount\PLy\newdimen\PLyy\newdimen\PLX
\newbox\PLdot \setbox\PLdot\hbox{{\tiny.}} \def\scl{.07} 
\def\PLot#1{\PLx`#1\advance\PLx-42\PLy\PLx\PLv\PLx\divide\PLy9\PLw\PLy\multiply
\PLw9\advance\PLx-\PLw\advance\PLx-4\PLy-\PLy\advance\PLy4\PLX=\the\PLx pt
\advance\PLyy\the\PLy pt\wd\PLdot=\scl\PLX\raise\scl\PLyy\copy\PLdot}
\def\draw#1{\ifx#1\end\let\next=\relax\else\PLot#1\let\next=\draw\fi\next}

\newcount\smallPLv\newcount\smallPLw\newcount\smallPLx\newcount\smallPLy\newdimen\smallPLyy\newdimen\smallPLX\newbox\smallPLdot \setbox\smallPLdot\hbox{{\tiny.}} \def\smallscl{.062} \def\smallPLot#1{\smallPLx`#1\advance\smallPLx-42\smallPLy\smallPLx\smallPLv\smallPLx\divide\smallPLy9\smallPLw\smallPLy\multiply \smallPLw9\advance\smallPLx-\smallPLw\advance\smallPLx-4\smallPLy-\smallPLy\advance\smallPLy4\smallPLX=\the\smallPLx pt \advance\smallPLyy\the\smallPLy pt\wd\smallPLdot=\smallscl\smallPLX\raise\smallscl\smallPLyy\copy\smallPLdot} \def\smalldraw#1{\ifx#1\end\let\next=\relax\else\smallPLot#1\let\next=\smalldraw\fi\next}

\newlength{\parrdp}\newlength{\parrht}
\newcommand{\parr}{\raisebox{-\parrdp}{\raisebox{\parrht}{\rotatebox{180}{$\&$}}}}
\newlength{\smallparrdp}\newlength{\smallparrht}
\newcommand{\smallparr}
{\mkern1mu\raisebox{-.1\smallparrdp}{\raisebox{\smallparrht}{\rotatebox{180}{\small$\&$}}}\mkern1mu}
\newlength{\footnoteparrdp}\newlength{\footnoteparrht}
\newcommand{\footnoteparr}
{\mkern1mu\raisebox{-\footnoteparrdp}{\raisebox{\footnoteparrht}{\rotatebox{180}{\footnotesize$\&$}}}\mkern1mu}
\newlength{\scriptparrdp}\newlength{\scriptparrht}
\newlength{\tinyparrdp}\newlength{\tinyparrht}
\newcommand{\C}{\mathbb{C}}

\newcommand{\N}{\mathcal{N}}

\newcommand{\tensor}{\otimes}

\newcommand{\dual}[1]{#1\perp}

\newcommand{\id}{\mathsf{id}}

\renewcommand{\perp}{^\bot}
\newcommand{\perpp}{{}^\bot{}^\bot}

\newcommand{\rulelabel}[1]{\mathsf{#1}}
\newcommand{\tensorlabel}{\otimes}
\newcommand{\parlabel}{\parr}
\newcommand{\onelabel}{1}
\newcommand{\botlabel}{\bot}
\newcommand{\axlabel}{\rulelabel{ax}}
\newcommand{\cutlabel}{\rulelabel{\mkern1mu cut}}
\newcommand{\smallparlabel}{\smallparr}

\psset{radius=2pt,xunit=.5cm,yunit=1,arrowscale=1.7,arrowlength=.7,labelsep=3pt,nodesep=.5ex}

\newcommand{\linevecpos}[3]{\ncline[ArrowInside=->,ArrowInsidePos=#3]{#1}{#2}}
\newcommand{\linevec}[2]{\linevecpos{#1}{#2}{.54}}
\newcommand{\vecanglesposheight}[6]{\nccurve[ncurv=#6,ArrowInside=->,ArrowInsidePos=#5,angleA=#3,angleB=#4]{#1}{#2}}
\newcommand{\vecanglespos}[5]{\vecanglesposheight{#1}{#2}{#3}{#4}{#5}{1}}
\newcommand{\vecanglesheight}[5]{\vecanglesposheight{#1}{#2}{#3}{#4}{.54}{#5}}
\newcommand{\vecangles}[4]{\vecanglespos{#1}{#2}{#3}{#4}{.54}}
\newcommand{\uvec}[2]{\vecangles{#1}{#2}{90}{-90}}
\newcommand{\dvec}[2]{\vecangles{#1}{#2}{-90}{90}}

\newcommand{\dloopvecanglesheight}[5]{\nccurve[ncurv=#5,angleA=#3,angleB=#4]{#1}{#2}
\lput{:0}{\psline{->}(.05,.005)(.1,.005)}}
\newcommand{\dloopvecheight}[3]{\dloopvecanglesheight{#1}{#2}{-75}{-105}{#3}}
\newcommand{\dloopvec}[2]{\dloopvecheight{#1}{#2}{1.3}}

\newcommand{\uloopvecanglesheight}[5]{\nccurve[ncurv=#5,angleA=#3,angleB=#4]{#1}{#2}
\lput{:0}{\psline{->}(.05,-.005)(.1,-.005)}}
\newcommand{\uloopvecheight}[3]{\uloopvecanglesheight{#1}{#2}{75}{105}{#3}}
\newcommand{\uloopvecleft}[2]{\nccurve[ncurv=1.3,angleA=105,angleB=75]{#1}{#2}
\lput{:0}{\psline{->}(.05,-.005)(.1,-.005)}}

\newcommand{\goi}{\mathbf{GoI}}

\newcommand{\one}[1]{\rnode{#1}{1}}
\newcommand{\nbot}[1]{\rnode{#1}{\bot}}

\newcommand{\pfunc}{\mathsf{Setp}}
\newcommand{\finsetp}{\mathsf{fSetp}}
\newcommand{\goipfunc}{\goi(\pfunc)}
\newcommand{\goifinsetp}{\goi(\finsetp)}

\newcommand{\bgap}{\rule{9mm}{0mm}}
\newcommand{\sgap}{\rule{9mm}{0mm}}
\newcommand{\negb}[1]{\cnode{2pt}{#1}}
\newcommand{\posb}[1]{\cnode*{2pt}{#1}}
\newcommand{\gnegb}[1]{\bgap\negb{#1}}
\newcommand{\gposb}[1]{\bgap\posb{#1}}
\newcommand{\posrestr}{^+}
\newcommand{\negrestr}{^-}

\newcommand{\cutpair}[1]{{\psset{nodesep=1pt}\rnode{a}{#1}\;\,\rnode{b}{#1}\perp%
\nccurve[ncurv=1,angleA=-60,angleB=-120]{a}{b}}}
\newcommand{\reduc}[1]{|#1|}
\newcommand{\Matching}{\textsc{Matching}}
\newcommand{\Switching}{\textsc{Switching}}
\newcommand{\parseone}{3.1ex}
\newcommand{\parsetwo}{7.6ex}
\newcommand{\parsethree}{12.1ex}
\newcommand{\parsenodeparents}[7]{\ncline[linestyle=none]{#1}{#2}\nbput[labelsep=#3]{\rnode{#4}{#5}}%
\ncline{#6}{#4}\ncline{#7}{#4}}
\newcommand{\parsenode}[5]{\parsenodeparents{#1}{#2}{#3}{#4}{#5}{#1}{#2}}
\newcommand{\parsenodeone}[4]{\parsenode{#1}{#2}{\parseone}{#3}{#4}}

\newcommand{\parsenodetwoleft}[5]{\parsenodeparents{#1}{#2}{\parsetwo}{#3}{#4}{#1}{#5}}
\newcommand{\parsenodetworight}[5]{\parsenodeparents{#1}{#2}{\parsetwo}{#3}{#4}{#5}{#2}}
\newcommand{\parsenodethreeleft}[5]{\parsenodeparents{#1}{#2}{\parsethree}{#3}{#4}{#1}{#5}}
\newcommand{\parsenodethreeright}[5]{\parsenodeparents{#1}{#2}{\parsethree}{#3}{#4}{#5}{#2}}
\newcommand{\sync}{\Bumpeq}
\newcommand{\lam}{\mathbf{Lam}}
\newcommand{\synccomp}{\mkern1mu\textbf{;}\mkern1mu}
\newcommand{\im}{\,\text{im}\,}
\newcommand{\dom}{\,\text{dom}\,}

\title{\vspace*{-5.5ex}\Large\bf 
Simple multiplicative proof nets with units
}%
\author{\\[-3.7ex]\normalsize\sc 
Dominic J.~D.~Hughes
\\[-.2ex]
\normalsize Stanford University\\
\small March 2005\thanks{Submitted for publication, March 5, 2005.  Comments welcome.}}
\date{}
\begin{document}
\maketitle\vspace*{-4ex}

\begin{abstract}
This paper presents a simple notion of proof net for multiplicative
linear logic with units.  Cut elimination is direct and strongly
normalising, in contrast to previous approaches which resorted to
moving jumps (attachments) of par units during normalisation.
Composition in the resulting category of proof nets is simply path
composition: all of the dynamics happens in $\goipfunc$, the
geometry-of-interaction construction applied to the category of sets
and partial functions.

\vspace{1ex}\noindent\textbf{Keywords: }multiplicative linear logic, units, proof nets, geometry of interaction.

\vspace{1ex}\noindent\textbf{AMS subject classification: }3B47
Substructural logics, 03F52 Linear logic.
\end{abstract}

\section{Introduction}

Here is a passage from Girard's \emph{Proof Nets: the Parallel Syntax
for Proof Theory} \cite[\S A.2]{Gir96}\footnote{Similar remarks are
in the earlier \emph{Linear Logic: A Survey} \cite[\S3.6]{Gir93}.}:
\begin{quote}\label{Gir-quote}\small
There are two multiplicative neutrals, $1$ and $\bot$, and two rules,
the axiom $\vdash 1$ and the weakening rule: from $\vdash\Gamma$,
deduce $\vdash\Gamma,\bot$.  Both rules are handled by means of links
with one conclusion and no premise; however, $\bot$-links are treated
like $0$-ary $?$-links, \ie, they must be given a default jump.
Sequentialisation is immediate.

At first sight, cut elimination is unproblematic: replace a cut
between conclusions $1$ and $\bot$ of zero-ary links with\ldots
nothing.  But we notice a new problem, namely that a cut formula $A$
can be the default jump of a $\bot$-link $L$, and we must therefore
propose another jump for $L$.  Usually one of the premises of the link
with conclusion $A$ works (or the jump of $L'$ if $A$ is the
conclusion of a $\bot$-link) works.  Worse, this new jump is by no
means natural (if $A$ is $B\tensor C$, the new jump can either be $B$
or $C$), which is quite unpleasant.  As far as we know, the only
solution consists in declaring that jumps are not part of the
proof-net, but rather some control structure.  It is then enough to
show that at least one choice of default jump is possible.  This is
not a very elegant solution: we are indeed working with equivalence
classes of proof nets and if we want to be rigorous we shall have to
endlessly check that such and such operation does not depend on the
choice of default jumps.  
\end{quote}

This paper presents a very simple solution: define a multiplicative
proof net with units (neutrals) as a function from negative to
positive formula leaves, satisfying the usual correctness criterion
\cite{Gir87,DR89}.  Cut elimination on binary connectives is then
trivial (as usual in the unit-free setting), and we have a direct
strong normalisation by standard path composition: all of the dynamics
happens in $\goipfunc$, the geometry-of-interaction or feedback
construction
\cite{Gir89,JSV96,Abr96} applied to the category
of sets and partial functions.

The novelty here is not the directed edges between negative and
positive leaves, an idea which goes back to the origins of linear
logic \cite{Gir87} and Kelly-MacLane graphs \cite{KM71}.  The key
contribution is the simply defined, strongly normalising cut
elimination, over $\goipfunc$.

\parag{The nets} Here is a simple example of a cut-free proof net on a four-formula sequent:
$$
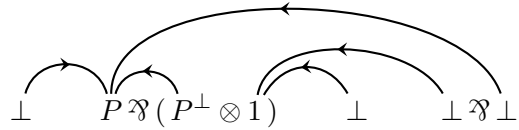
\begin{psmatrix}\label{example}
\nbot{b1}\sgap\rnode{P1}{P}\,\parr\,(\,\rnode{P2}{P}\perp\tensor\one{i}\,)\sgap\nbot{b2}\sgap\,\nbot{b3}\,\parr\,\nbot{b4}
\rule{0pt}{1.5cm}
\uloopvecheight{b1}{P1}{1}
\uloopvecanglesheight{P2}{P1}{105}{75}{1}
\vecanglesheight{b2}{i}{120}{60}{1}
\vecanglesheight{b3}{i}{120}{80}{.7}
\vecanglesheight{b4}{P1}{110}{90}{.6}
\end{psmatrix}$$
The graph of the function from negative to positive leaves is shown
by the directed edges.
Note that all four switchings are trees.  This is easier to see if we
show the parse trees:
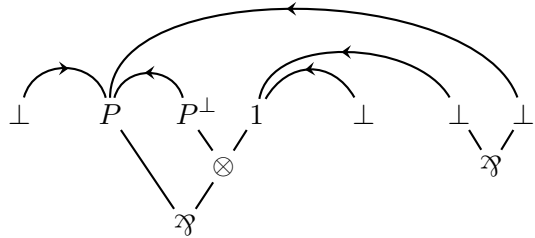
\begin{center}\begin{math}
\nbot{b1}\sgap\rnode{P1}{P}\,\;\;\;\;\;\,\,\rnode{P2}{P}\perp\;\;\;\;\one{i}\;\;\,\sgap\nbot{b2}\sgap\,\nbot{b3}\,\;\;\;\;\,\nbot{b4}\,
\rule{0pt}{1.5cm}\vspace*{7ex}
\uloopvecheight{b1}{P1}{1}
\uloopvecanglesheight{P2}{P1}{105}{75}{1}
\vecanglesheight{b2}{i}{120}{60}{1}
\vecanglesheight{b3}{i}{120}{80}{.7}
\vecanglesheight{b4}{P1}{110}{90}{.6}
\parsenodeone{P2}{i}{tensor}{\tensor}
\parsenodetwoleft{P1}{i}{par}{\parr}{tensor}
\parsenodeone{b3}{b4}{par2}{\parr}
\end{math}\end{center}
As with the unit-free case \cite{Gue99,MO00}, correctness can be
checked in linear time (see Section~\ref{linear}).
\parag{GoI dynamics}
MLL formulas and proof nets form a category with a morphism $A\to B$ a
cut-free proof net on $\vdash A\perp,B$.
For example,
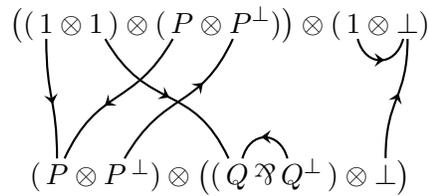
\begin{center}\label{firstexample}
\begin{psmatrix}
$\big((\,\one{u}\tensor\one{uu}\,)\tensor(\,\rnode{a}{P}\tensor\rnode{aa}{P}\perp)\big)
\tensor (\,\one{uuu}\tensor \nbot{uuuu})$ \\
$(\,\rnode{a1}{P}\tensor\rnode{aa1}{P}\,\perp) \tensor \big((\,\rnode{b1}{Q}\,\parr\,\rnode{bb1}{Q}\perp\,)
\tensor \nbot{u1}\big)$
\dvec{u}{a1}
\vecanglespos{uu}{b1}{-60}{115}{.5}
\vecanglespos{a}{a1}{-120}{65}{.59}
\vecanglespos{aa1}{aa}{60}{-115}{.67}
\uloopvecleft{bb1}{b1}
\dloopvec{uuu}{uuuu}
\uvec{u1}{uuuu}
\end{psmatrix}
\end{center}
is a morphism from the upper formula to the lower formula.  (We
suppress the negation on the input/upper formula, flipping polarity, so
tensors are switched in the input.)
The underlying $\goipfunc$ morphism is:
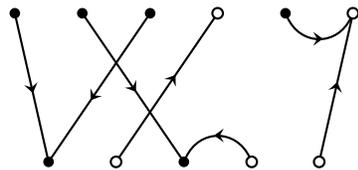
\begin{center}\vspace*{-3.6ex}
\begin{math}
\psset{nodesep=-.2mm}
\begin{psmatrix}
\posb{u}\bgap\posb{uu}\bgap\posb{a}\bgap\negb{aa}\bgap\posb{uuu}\bgap\negb{uuuu} \\
\posb{a1}\bgap\negb{aa1}\bgap\posb{b1}\bgap\negb{bb1}\bgap\negb{u1}
\linevec{u}{a1}
\linevecpos{uu}{b1}{.53}
\linevecpos{a}{a1}{.59}
\linevecpos{aa1}{aa}{.59}
\uloopvecanglesheight{bb1}{b1}{115}{65}{1}
\dloopvecanglesheight{uuu}{uuuu}{-65}{-115}{1}
\linevec{u1}{uuuu}
\end{psmatrix}
\end{math}
\end{center}
An object of $\goipfunc$ is a signed set $S$, whose elements we shall
call
\emph{leaves}, and a morphism $S\to T$ is a partial
function from 
negative leaves to positive leaves (polarity flipped on the input side).
Composition is standard path composition, \eg
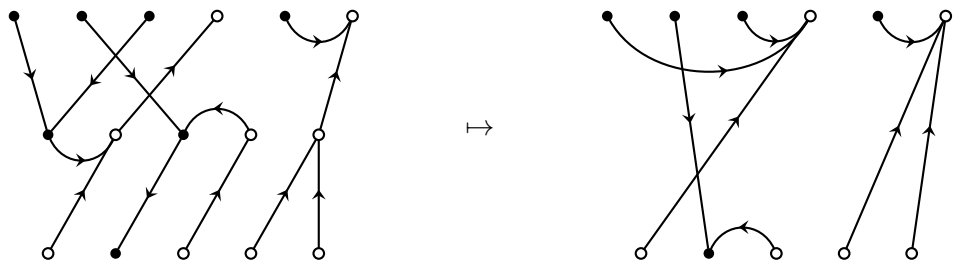
\begin{center}\label{goi-comp}
\psset{rowsep=1.1cm,nodesep=-.2mm}\begin{psmatrix}
$\posb{u}\bgap\posb{uu}\bgap\posb{a}\bgap\negb{aa}\bgap\posb{uuu}\bgap\negb{uuuu}$ 
&&
$\posb{u'}\bgap\posb{uu'}\bgap\posb{a'}\bgap\negb{aa'}\bgap\posb{uuu'}\bgap\negb{uuuu'}$
\\
$\posb{a1}\bgap\negb{aa1}\bgap\posb{b1}\bgap\negb{bb1}\bgap\negb{u1}$
&
$\mapsto$ \\
$\negb{u2}\gposb{b2}\gnegb{bb2}\gnegb{uu2}\gnegb{uuu2}$
&&
$\negb{u2'}\gposb{b2'}\gnegb{bb2'}\gnegb{uu2'}\gnegb{uuu2'}$
\linevec{u}{a1}
\linevecpos{uu}{b1}{.53}
\linevecpos{a}{a1}{.59}
\linevecpos{aa1}{aa}{.59}
\uloopvecanglesheight{bb1}{b1}{115}{65}{1}
\dloopvecanglesheight{uuu}{uuuu}{-65}{-115}{1}
\linevec{u1}{uuuu}
\dloopvecanglesheight{a1}{aa1}{-65}{-115}{1}
\linevec{u2}{aa1}
\linevec{b1}{b2}
\linevec{bb2}{bb1}
\linevec{uu2}{u1}
\linevec{uuu2}{u1}
\vecanglesposheight{u'}{aa'}{-60}{-120}{.58}{.8}
\linevecpos{uu'}{b2'}{.45}
\dloopvecanglesheight{a'}{aa'}{-65}{-120}{1}
\dloopvecanglesheight{uuu'}{uuuu'}{-65}{-115}{1}
\linevecpos{u2'}{aa'}{.58}
\uloopvecanglesheight{bb2'}{b2'}{115}{65}{1}
\linevec{uu2'}{uuuu'}
\linevec{uuu2'}{uuuu'}
\end{psmatrix}
\end{center}
which provides composition (turbo cut elimination) in the category of proof nets, \eg
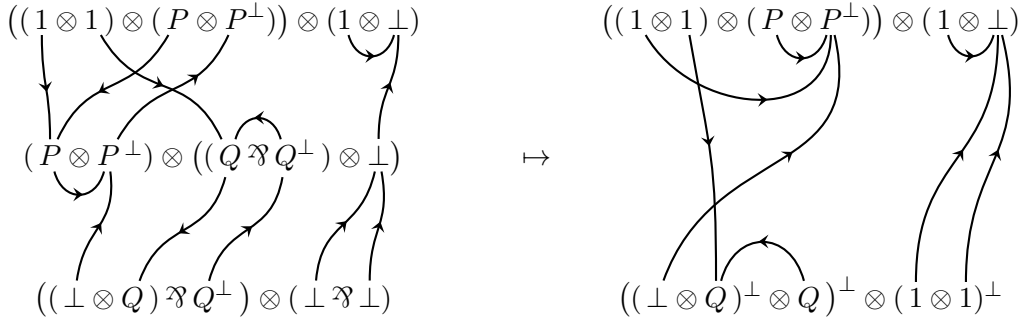
\begin{center}\label{firstexamplecomp}
\psset{rowsep=1.3cm,colsep=1cm}\begin{psmatrix}
$\big((\,\one{u}\tensor\one{uu}\,)\tensor(\,\rnode{a}{P}\tensor\rnode{aa}{P}\perp)\big)
\tensor (\one{uuu}\tensor \nbot{uuuu})$
&&
$\big((\,\one{u'}\tensor\one{uu'}\,)\tensor(\,\rnode{a'}{P}\tensor\rnode{aa'}{P}\perp)\big)
\tensor (\,\one{uuu'}\tensor \nbot{uuuu'})$
\\
$(\,\rnode{a1}{P}\tensor\rnode{aa1}{P}\,\perp) \tensor \big((\,\rnode{b1}{Q}\,\parr\,\rnode{bb1}{Q}\perp\,)
\tensor \nbot{u1}\big)$
&
\hspace{2ex}$\mapsto$\hspace{-2ex} \\
$\big((\,\nbot{u2}\tensor \rnode{b2}{Q}\,)\,\parr\, \rnode{bb2}{Q}\perp\,\big)
\tensor (\,\nbot{uu2} \,\parr\, \nbot{uuu2})$
&&
$\big((\,\nbot{u2'}\tensor \rnode{b2'}{Q}\,)\perp\tensor \rnode{bb2'}{Q}\,\big)\perp
\tensor (\,\one{uu2'} \tensor \one{uuu2'})\perp$
\dvec{u}{a1}
\vecanglespos{uu}{b1}{-60}{115}{.5}
\vecanglespos{a}{a1}{-120}{65}{.59}
\vecanglespos{aa1}{aa}{60}{-115}{.67}
\uloopvecleft{bb1}{b1}
\dloopvec{uuu}{uuuu}
\uvec{u1}{uuuu}
\dloopvec{a1}{aa1}
\vecangles{u2}{aa1}{82}{-83}
\vecangles{b1}{b2}{-105}{70}
\vecangles{bb2}{bb1}{75}{-100}
\vecangles{uu2}{u1}{85}{-103}
\vecangles{uuu2}{u1}{90}{-80}
\vecanglesposheight{u'}{aa'}{-60}{-90}{.58}{1}
\vecanglespos{uu'}{b2'}{-80}{90}{.45}
\dloopvec{a'}{aa'}
\dloopvec{uuu'}{uuuu'}
\vecanglespos{u2'}{aa'}{70}{-77}{.58}
\uloopvecleft{bb2'}{b2'}
\vecangles{uu2'}{uuuu'}{90}{-90}
\vecangles{uuu2'}{uuuu'}{90}{-76}
\end{psmatrix}
\end{center}
is the path composition of the previous $\goi$ diagram.
This provides a simple solution to the problems articulated by
Girard above.

\parag{Laminated GoI composition for MALL nets} Section~\ref{mall}
continues the $\goi$ theme, and shows how composition (turbo cut
elimination) of MALL proof nets \cite{HG03,HG05} can be viewed as
occurring in a `laminated' variant of $\goipfunc$.

\parag{Related work}  Proof nets with units are in \cite{BCST96} and \cite{LS04}.  
Neither solves the problems in Girard's quote: each suffers from the
need to move $\bot$-jumps during elimination, so
one is lumbered once again with equivalence classes.  The cut-free
one-sided MLL proof nets
in \cite{BCST96} are the cut-free proof nets described in Girard's
quote in a circuit/wire notation, with an additional ordering on
$\bot$-jumps:
see Section~\ref{Girard-BCST}.  
The paper \cite{LS04} defines a cut-free proof net on a sequent
$\vdash\Gamma$ as a separate MLL formula $\Theta$ whose leaves from
left-to-right are a permutation of those of $\Gamma$.
The $\bot$-jumps and axiom links are thus enveloped in an additional
syntactic layer $\Theta$: see Section~\ref{LS}.
The proof nets of \cite{MO03} for intuitionistic multiplicative linear
logic with units (based on essential nets \cite{Lam94}) involve
directed edges.

Work in progress quotients the nets presented in this paper by
Trimble's \emph{empire rewiring} \cite{Tri94}, which permits a
$\bot$-jump target to move so long as correctness is not broken, to
construct free star-autonomous categories for full coherence (\cf\
\cite{BCST96,KO99,MO03,LS04}).

\parag{Acknowledgement} Thanks to Robin Houston for feedback.

\section{Notation}

By MLL we mean multiplicative linear logic with units \cite{Gir87}.
Formulas are built from literals (propositional variables $P,Q,\ldots$
and their duals $\dual{P}$, $\dual{Q},\ldots$) and
units/constants/neutrals $1$ and $\bot$ by the binary connectives
\defn{tensor}~$\tensor$ and \defn{par}~$\parr$.
Negation\label{negation} $(-)\perp$ extends to arbitrary formulas with
$P\perpp\mkern-2mu=\mkern-1mu P$ on propositional variables,
$\;\bot\perp\mkern-2mu=\mkern-2mu 1\,$, $\;1\perp\mkern-2mu=\mkern-2mu\bot\,$, and de
Morgan duality $(A\tensor B)\perp=A\perp\parr B\perp$ and $(A\parr
B)\perp=A\perp\tensor B\perp$.  An \defn{atom} is a literal or unit.
We identify a formula with its parse tree: a tree labelled with
atoms at the leaves and connectives at internal vertices.  
A \defn{sequent} is a non-empty disjoint union of formulas.  Thus a
sequent is a particular kind of labelled forest.  We write comma for
disjoint union.
Sequents are proved using the following rules:
{$$
\begin{prooftree}\thickness=.08em
\strut
\justifies
\,P, P\perp
\using\axlabel
\end{prooftree}
\hspace{4ex}
\begin{prooftree}\thickness=.08em
\Gamma,A\;\;\;\;\;A\perp,\Delta
\justifies
\Gamma,\Delta
\using\cutlabel
\end{prooftree}
\hspace{4ex}
\begin{prooftree}\thickness=.08em
\rule{0pt}{1.4ex}
\justifies
\;\;1\;\using\onelabel
\end{prooftree}
\hspace{4ex}
\begin{prooftree}\thickness=.08em
\Gamma
\justifies
\;\Gamma,\,\bot\using\botlabel
\end{prooftree}
\hspace{4ex}
\begin{prooftree}\thickness=.08em
\Gamma,A\;\;\;\;\;B,\Delta
\justifies
\Gamma,A\tensor B,\Delta
\using\tensorlabel
\end{prooftree}
\hspace{4ex}
\begin{prooftree}\thickness=.08em
\Gamma,\,A,\,B
\justifies
\Gamma,\,A\smallparr B
\using\smallparlabel
\end{prooftree}
$$}%
Here, and throughout this document, $P\mkern-2mu,Q,\ldots$ range over
propositional variables, $A,B,\ldots$ over formulas, and
$\Gamma\mkern-3mu,\Delta,\ldots$ over (possibly empty) disjoint unions of
formulas.
Without loss of generality 
we restrict the axiom rule to literals \cite{Gir87}.
The propositional variables $P,Q,\ldots$ and the unit $1$ are
\defn{positive}, and their duals $P\perp,Q\perp,\ldots$ and $\bot$ are \defn{negative}.  A leaf 
of a formula is positive/negative according to its label.
A \defn{cut pair} $\cutpair{A}$ is a disjoint union of complementary
formulas $A$ and $A\perp$ together with an undirected edge, a
\defn{cut}, between their roots.
A \defn{cut sequent} is a disjoint union of a sequent and zero or more
cut pairs.
A \defn{switching} of a cut sequent is any subgraph obtained by
deleting one of the two argument edges of each $\parr$ (see
\cite{DR89}).
By an \defn{old proof net} we mean a proof net for MLL with units
as in Girard's quote in the Introduction; see
\cite{Dan90,Reg92,GSS92,Gir93,Gir96} for history and development. (An example of an old proof net is drawn in 
the next section.)

\section{Proof nets}

A \defn{leaf function} on a cut sequent is a function from its negative leaves to its positive leaves.
A \defn{proof net} on a cut sequent $\Gamma$ is a 
leaf function $f$ on $\Gamma$ satisfying:
\begin{itemize}\vspace*{-.8ex}
\item \Matching. For any propositional variable $P$, the restriction of $f$ to $P$-labelled leaves is a bijection 
between the $P$-labelled leaves of $\Gamma$ and the $P\perp$-labelled
leaves of $\Gamma$.
\item \Switching. For any switching $\Gamma'$ of $\Gamma$, the undirected graph obtained by adding the 
edges of $f$ to $\Gamma'$ is a tree (acyclic and connected).
\end{itemize}
See page~\pageref{example} for an example.
This definition amounts to a restricted case of an old proof net:
restrict $\bot$-jumps to target positive leaves and reject unit axiom
links (use $\bot\to 1$ jumps instead).  In addition, we orient all
axiom links from negative to positive.
Stating this the other way round, the above definition relaxes to the
old definition thus: (a) on $\bot$-labelled leaves allow $f$ to target
any vertex (equivalently subformula) of $\Gamma$, not just a positive
leaf, (b) distinguish between two kinds of edges from $\bot$ to $1$
(jump \emph{versus} axiom link), and (c) draw axiom links unoriented.
Here is an example of an old proof net:
$$
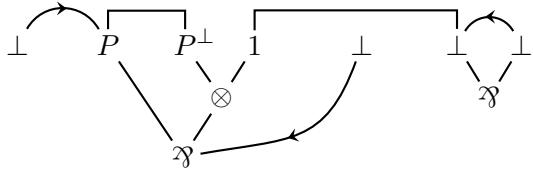
\begin{psmatrix}\label{oldnet}
\nbot{b1}\sgap\rnode{P1}{P}\,\;\;\;\;\;\,\,\rnode{P2}{P}\perp\;\;\;\;\one{i}\;\;\,\sgap\nbot{b2}\sgap\,\nbot{b3}\,\;\;\;\;\,\nbot{b4}\,
\rule{0pt}{.6cm} \\
\parsenodeone{P2}{i}{tensor}{\tensor}
\parsenodetwoleft{P1}{i}{par}{\parr}{tensor}
\parsenodeone{b3}{b4}{par2}{\parr}
\vecanglesheight{b2}{par}{-110}{10}{1}
\ncbar[angleA=90,angleB=90,arm=7pt]{P2}{P1}
\ncbar[angleA=90,angleB=90,arm=7pt]{i}{b3}
\uloopvecanglesheight{b1}{P1}{60}{120}{.9}
\psset{nodesep=1.5pt}
\uloopvecanglesheight{b4}{b3}{120}{60}{.9}
\end{psmatrix}\vspace*{-3ex}$$
which in original proof net notation is:
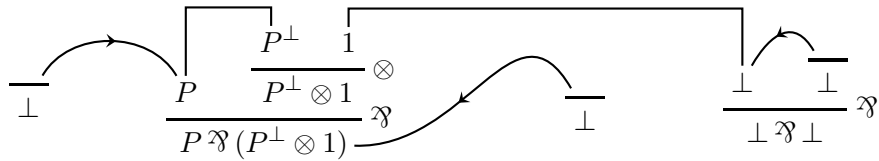
\begin{center}\label{oldnetold}\vspace*{2ex}
\begin{prooftree}\thickness=.08em
\justifies\rnode{b1}{\;\bot}\;
\end{prooftree}
\hspace{8ex}
\begin{prooftree}\thickness=.08em
\;\rnode{P1}{P}\;\;\;\;\;
 \[
  \;\rnode{P2}{P}\perp\;\;\;\;\;\rnode{i}{1}\;
  \justifies
  \;P\perp\tensor 1\; \using \tensorlabel
 \]
\justifies
\;\rnode{par}{P\,\parr\,(P\perp\tensor 1)}\;
\using \parlabel
\end{prooftree}
\hspace{12ex}
\begin{prooftree}\thickness=.08em
\strut\justifies\;\rnode{b2}{\bot\;}
\end{prooftree}
\hspace{8ex}
\begin{prooftree}\thickness=.08em
\;\nbot{b3}\;\;\;\;\;\[ \justifies\;\rnode{b4}{\bot\;} \]
\justifies
\bot\,\parr\,\bot \using \parlabel
\end{prooftree}
\ncbar[angleA=90,angleB=90,arm=7pt]{P2}{P1}
\ncbar[angleA=90,angleB=90,arm=7pt]{i}{b3}
\psset{nodesepA=8pt}
\vecanglesheight{b2}{par}{120}{-2}{1}
\vecanglesheight{b1}{P1}{60}{120}{.8}
\vecanglesheight{b4}{b3}{120}{60}{1.2}
\end{center}
Axiom links are shown as three-segment straight edges, and
jumps from $\bot$-links
$\overline{\,\bot\,}$ are shown curved and directed.

Translation from a proof to a proof net is as usual, with a
$\bot$-jump added at each $\bot$-rule, but now with choice of target
restricted to positive atoms only.  Note that well-definedness relies
on the observation that every provable MLL sequent contains a positive
atom.  The translation becomes deterministic upon marking a positive
leaf in the conclusion of every $\bot$-rule.
For example, each of the following marked proofs translates
(deterministically) into the proof net on page~\pageref{example}:
\begin{center}
\begin{prooftree}\thickness=.08em
\[
  \[
   \[
    \[
     \[
      \;\;\:
      \[
       \justifies
       \;P\:,\:P\perp\;
       \using \axlabel
      \]
      \justifies
      \;\bot,\,P\;,\;P\perp\;
      \using \botlabel
     \]
     \;\;
     \[
      \justifies
      \; 1\,,\,\bot\;
      \using \axlabel
     \]
     \justifies
     \;\bot,\,P\;,\;\,P\perp\mkern-2mu\tensor\mkern-2mu 1\;\,,\,\bot\;
     \using \tensorlabel
    \]
    \mbox{}\hspace{-2.7ex}
    \justifies
    \;\bot,\,P\;,\;\,P\perp\mkern-2mu\tensor\mkern-2mu\underline 1\;\,,\,\bot,\,\bot\;
    \using \botlabel
   \]
   \justifies
   \;\bot,\,P\mkern1mu\parr(P\perp\mkern-2mu\tensor\mkern-2mu 1),\,\bot,\,\bot\;
   \using \parlabel
  \]
  \mbox{}
  \justifies
  \;\bot,\,\underline{P}\mkern1mu\parr(P\perp\mkern-2mu\tensor\mkern-2mu 1),\,\bot,\,\bot\,,\,\bot\;
  \using \botlabel
\]
\justifies 
\;\bot,\,P\mkern1mu\parr(P\perp\mkern-2mu\tensor\mkern-2mu 1),\,\bot,\,\bot\parr\bot\;
\using\parlabel
\end{prooftree}
\hspace{15ex}
\begin{prooftree}\thickness=.08em
\;
\[
 \[
  \[
   \[
    \[
     \justifies
     \;P\;,\;P\perp\;
     \using \axlabel
    \]
    \;\;
    \[
     \[
       \[
       \justifies
       \;1\;
       \using \onelabel
       \]
       \:\mkern1mu
      \justifies
      \; 1\,,\,\bot\;
      \using \botlabel
     \]
     \mbox{}\hspace{-.8ex}
     \justifies
     \; 1\,,\,\bot\,,\,\bot
     \using \botlabel
    \]
    \justifies
    \;P\;,\;\,P\perp\mkern-2mu\tensor\mkern-2mu 1\;\,,\,\bot,\,\bot\;
    \using \tensorlabel
   \]
   \mbox{}\hspace{-2.4ex}
   \justifies
   \;\underline P\;,\;\,P\perp\mkern-2mu\tensor\mkern-2mu 1\;\,,\,\bot,\,\bot\,,\,\bot\;
   \using \botlabel
  \]
  \justifies
  \; P\;,\;\,P\perp\mkern-2mu\tensor\mkern-2mu 1\;\,,\,\bot,\,\bot\parr\bot\;
  \using \parlabel
 \]
 \justifies
 \;P\mkern1mu\parr(P\perp\mkern-2mu\tensor\mkern-2mu 1),\,\bot,\,\bot\parr\bot\;
 \using\parlabel
\]
\justifies 
\;\bot,\,\underline{P}\mkern1mu\parr(P\perp\mkern-2mu\tensor\mkern-2mu 1),\,\bot,\,\bot\parr\bot\;
\using\botlabel
\end{prooftree}
\end{center}
Marks are shown by underlining; when a sequent has just one positive
atom, we leave the mark implicit.  (Downward tracking of $\bot$'s is
vertical, except through the tensor rule.)

\begin{theorem}[Sequentialisation]
Every proof net is a translation of a proof.
\end{theorem}
This is simply a restriction of the theorem for old proof nets.
Correctness is verifiable in linear time (a simple corollary of the
unit-free case \cite{Gue99,MO00}): see Section~\ref{linear}.
\section{Cut elimination}

Let $f$ be a proof net on the cut sequent $\Gamma,\cutpair{A}$.
The result $f'$ of \defn{eliminating} the cut in $\cutpair{A}$
is:
\begin{itemize}
\item \emph{Atom.}  Suppose $A$ is an atom.  Without loss of generality, $A$ is positive.  Delete $\cutpair{A}$ 
and reset any $f$-edge to $A$ to target $f(A\perp)$
instead.
\item \emph{Compound.} Suppose $A$ is a compound formula.  
Without loss of generality $A=B\tensor C$ and $A\perp=B\perp\parr
C\perp$.  Replace $\cutpair{A}$ by $\cutpair{B},\cutpair{C}$.  The
leaves, and $f$, remain unchanged.
\end{itemize}
Schematically:
\begin{center}
\hspace{8ex}\psset{colsep=3cm,framesep=0pt,nodesep=0pt,rowsep=1cm}
\begin{psmatrix}
\Cnode(-3,.5){a}
\Cnode(-2.5,.9){b}
\Cnode(-1.7,1.2){c}
\circlenode[linestyle=none]{A}{$A$}\rule{0pt}{7.2ex}
\hspace{4ex}
\circlenode[linestyle=none]{AA}{$A\perp\hspace{-1ex}$}
{{\psset{nodesep=1pt}\nccurve[ncurv=1,angleA=-60,angleB=-120]{A}{AA}}}
\psset{nodesepB=-1pt}
\vecangles{a}{A}{10}{150}
\vecangles{b}{A}{0}{125}
\vecangles{c}{A}{-10}{100}
\hspace{4ex}
\Cnode*(1,.2){d}
\psset{nodesepB=-1pt,nodesepA=-1pt}
\vecangles{AA}{d}{60}{110}
&
\raisebox{1ex}{$
\pstree[treemode=U,nodesep=2pt,levelsep=20pt]
{\TR[name=t]{\tensor}}
  {\TR[name=b]{B}
   \TR[name=c]{C}}
\hspace{8ex}
\pstree[treemode=U,nodesep=2pt,levelsep=20pt]
{\TR[name=p]{\parr}}
  {\TR[name=bb]{B\perp\!\!}
   \TR[name=cc]{C\perp\!\!}}
$}
\nccurve[nodesepA=2pt,nodesepB=1pt,angleA=-40,angleB=-130]{t}{p}
\\
\rput{270}(-1,.2){$\mapsto$}\rput[l](-.5,.2){\small$A$ atomic}
&
\rput{270}(0,0){$\mapsto$}
\\
\Cnode(-3,.5){a}
\Cnode(-2.5,.9){b}
\Cnode(-1.7,1.2){c}
\hspace{14ex}
\psset{nodesepB=-1pt}
\Cnode*(1,.2){d}
\psset{nodesepB=-1pt,nodesepA=-1pt}
\vecangles{a}{d}{10}{160}
\vecangles{b}{d}{10}{145}
\vecangles{c}{d}{10}{130}
&
\raisebox{5ex}{$
\hspace{.5ex}\rnode{b}{B}
\hspace{4ex}
\rnode{c}{C}
\hspace{8ex}
\rnode{bb}{B}\perp
\hspace{4ex}
\rnode{cc}{C}\perp\hspace{-2ex}
\nccurve[nodesepA=1pt,nodesepB=2pt,angleA=-40,angleB=-130]{b}{bb}
\nccurve[nodesepA=1pt,nodesepB=2pt,angleA=-40,angleB=-130]{c}{cc}
$}
\end{psmatrix}
\end{center}
\begin{theorem}
Cut elimination is well-defined: eliminating a cut from a proof net yields a proof net.
\end{theorem}
\begin{proof}
The atomic case is trivial, since switchings and cycles correspond
before and after the elimination.  The compound case is the same as
the usual unit-free elimination \cite{Gir87,DR89,Gir93}.
\end{proof}
\begin{proposition}
Cut elimination is locally confluent.
\end{proposition}
\begin{proof}
The only non-trivial case is a pair of atomic eliminations.  This case
is clear from the following schematic involving two interacting atomic
cut redexes $\cutpair{A}$ and $\cutpair{B}$.
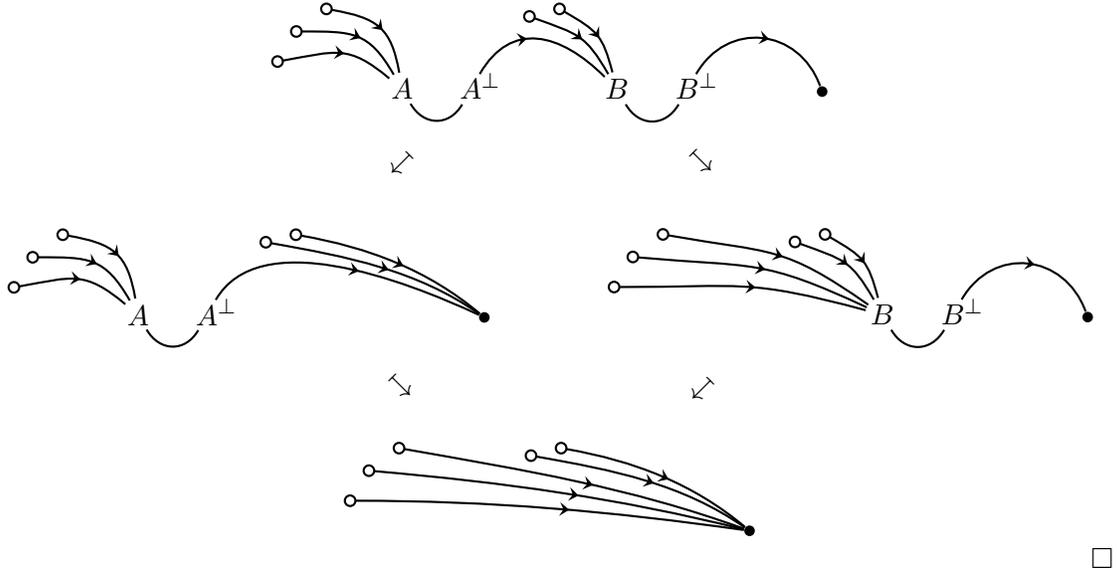
\begin{center}\psset{nodesep=-.2pt}
\newcommand{\abc}{\Cnode(-3,.5){a}
\Cnode(-2.5,.9){b}
\Cnode(-1.7,1.2){c}}
\newcommand{\bbcc}{\Cnode(-2,1.1){bb}
\Cnode(-1.2,1.2){cc}}
\newcommand{\abcA}{{\psset{nodesepB=2pt}\vecangles{a}{A}{10}{140}
\vecangles{b}{A}{0}{120}
\vecangles{c}{A}{-10}{100}}}
\newcommand{\bbccB}{{\psset{nodesepB=2pt}\vecangles{bb}{B}{-20}{120}
\vecangles{cc}{B}{-30}{105}}}
\newcommand{\maingap}{\hspace{8ex}}
\newcommand{\cutgap}{\hspace{3ex}}
\newcommand{\leftcut}{
$\circlenode[linestyle=none,framesep=-2pt]{A}{A}$
\cutgap
$\rnode{AA}{A}\perp$
\nccurve[ncurv=1,nodesep=2pt,angleA=-60,angleB=-120]{A}{AA}}
\newcommand{\rightcut}{
$\circlenode[linestyle=none,framesep=-2pt]{B}{B}$
\cutgap
$\rnode{BB}{B}\perp$
\nccurve[ncurv=1,nodesep=2pt,angleA=-60,angleB=-120]{B}{BB}}
\newcommand{\midjump}{{\psset{nodesep=2pt}\vecanglespos{AA}{B}{60}{135}{.4}}}
\newcommand{\abcB}{{\psset{nodesepB=2pt}\vecangles{a}{B}{0}{165}
\vecangles{b}{B}{-5}{155}
\vecangles{c}{B}{-10}{145}}}
\renewcommand{\C}{\maingap\Cnode*(0,.1){C}}
\newcommand{\rightjump}{{\psset{nodesepA=2pt}\vecangles{BB}{C}{60}{110}}}
\newcommand{\bbccC}{\vecanglesheight{bb}{C}{-10}{145}{.7}
\vecanglesheight{cc}{C}{-10}{140}{.7}}
\newcommand{\AAbbccC}{{\psset{nodesepA=2pt}\vecanglesheight{AA}{C}{60}{155}{.65}}\bbccC}
\newcommand{\abcC}{\vecanglesheight{a}{C}{0}{172}{.7}
\vecanglesheight{b}{C}{-5}{166}{.7}
\vecanglesheight{c}{C}{-10}{160}{.7}}
\newcommand{\finaledges}{\bbccC\abcC}
\begin{pspicture}(-3,.4)(0,7.4)
\rput(0,6){%
\abc
\leftcut
\abcA
\maingap
\bbcc
\rightcut
\midjump
\bbccB
\C
\rightjump
}
\rput(-8,3){%
\abc
\leftcut
\abcA
\maingap
\bbcc
\cutgap
\C
\AAbbccC
}
\rput(8,3){%
\abc
\cutgap
\maingap
\bbcc
\rightcut
\abcB
\bbccB
\C
\rightjump
}
\rput(0,0){%
\abc
\cutgap
\maingap
\bbcc
\cutgap
\C
\finaledges
}
\rput{-135}(-5.5,5){$\mapsto$}
\rput{-45}(2.5,5){$\mapsto$}
\rput{-135}(2.5,2){$\mapsto$}
\rput{-45}(-5.5,2){$\mapsto$}
\end{pspicture}
\end{center}
\end{proof}

\begin{theorem}
Cut elimination is strongly normalising.
\end{theorem}
\begin{proof}
It is locally confluent, and eliminating a cut reduces the number of
vertices of the cut sequent.
\end{proof}

\parag{Turbo cut elimination}

As with standard unit-free MLL proof nets, normalisation can be
completed in a single step.  For $l$ the $i\nth$ leaf of a formula $A$
in a cut pair $\cutpair{A}$, let $l\perp$ denote the $i\nth$ leaf of
$A\perp$.
The \defn{normal form} of a cut sequent $\Gamma$ is the sequent
$\reduc{\Gamma}$ obtained by deleting all cut pairs. Given a proof net
$f$ on $\Gamma$, its \defn{normal form} $\reduc{f}$ is the proof net
on $\reduc{\Gamma}$ obtained by replacing every set of edges $\langle
l_0,l_1\rangle,\langle l_1\perp,l_2\rangle,\langle
l_2\perp,l_3\rangle,\ldots,
\langle l_{n-1}\perp,l_n\rangle$ in $f$ in which only $l_0$
and $l_n$ occur in $\reduc{\Gamma}$ by the single edge $\langle
l_0,l_n\rangle$.  By a simple induction on the number of vertices of
cut sequents, $\reduc{f}$ is precisely the normal form of $f$ under
one-step cut elimination.  (In particular, this implies $\reduc{f}$ is
indeed a proof net.)

\section{GoI dynamics}

Cut elimination yields a category $\N$ of MLL proof nets.
Objects are MLL formulas, and a morphism $A\to B$ is a proof net on
the (cut-free) sequent $A\perp,B$ (\cf\ \cite{HG03,HG05}, for
example).  The composite of $f:A\to B$ and $g:B\to C$ is the normal
form of the proof net $f\cup g$ on $A\perp,\cutpair{B},C$.
Composition is associative because cut elimination is strongly
normalising.
The identity $A\to A$, a leaf function on $A\perp,A$, has an edge
between the $i\nth$ leaf of $A\perp$ and the $i\nth$ leaf of $A$ for
each $i$, oriented from negative to positive.

We generally draw $f:A\to B$ 
with $A$ above $B$, and suppress the negation on $A$.  For example,
the identity $\bot\tensor P\to \bot\tensor P$
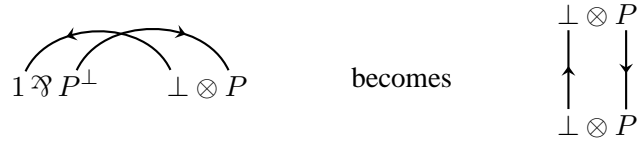
$$\begin{psmatrix}
\one{1}\,\parr\,\rnode{pp}{P}\perp\sgap\nbot{b}\tensor\rnode{p}{P}
\rule{0pt}{1cm}
\vecanglesposheight{b}{1}{115}{65}{.68}{.8}
\vecanglesposheight{pp}{p}{65}{115}{.68}{.8}
\end{psmatrix}
\hspace{8ex}\text{becomes}\hspace{8ex}
\raisebox{-3ex}{$\begin{psmatrix}[rowsep=1cm]
\nbot{b}\tensor\rnode{p}{P} \\
\nbot{b'}\tensor\rnode{p'}{P}
\uvec{b'}{b}\dvec{p}{p'}
\end{psmatrix}$}\hspace{5ex}
$$
Similarly, a composition such as
$$
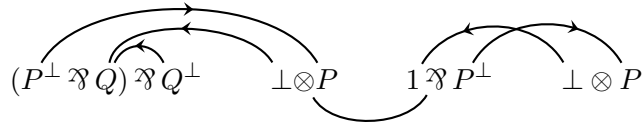
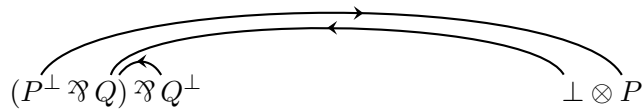
\begin{psmatrix}[rowsep=.8cm]
(\rnode{Pp}{P}\perp\,\parr\,\rnode{Q}{Q})\,\parr\,\rnode{Qp}{Q}\perp\sgap
\nbot{b}\rnode{tensor}{\tensor}\rnode{P}{P}\sgap\one{i}\,\rnode{par}{\parr}\,\rnode{Pp'}{P}\perp
\sgap\nbot{b'}\tensor \rnode{P'}{P}\rule{0pt}{1cm}\hspace{2ex}
\uloopvecanglesheight{Qp}{Q}{115}{66}{1}
\uloopvecanglesheight{b}{Q}{120}{80}{.6}
\vecanglesheight{Pp}{P}{60}{120}{.6}
\vecanglesposheight{b'}{i}{120}{60}{.68}{.8}
\vecanglesposheight{Pp'}{P'}{60}{120}{.68}{.8}
\nccurve[angleA=-60,angleB=-120]{tensor}{par}
\\
\rput{-90}{\mapsto}
\\
(\rnode{Pp}{P}\perp\,\parr\,\rnode{Q}{Q})\,\parr\,\rnode{Qp}{Q}\perp\sgap
\hspace{6ex}\sgap\hspace{6ex}
\sgap\nbot{b'}\tensor \rnode{P'}{P}\rule{0pt}{1cm}\hspace{2ex}
\uloopvecanglesheight{Qp}{Q}{120}{55}{1}
\uloopvecanglesheight{b'}{Q}{120}{75}{.3}
\vecanglesheight{Pp}{P'}{60}{120}{.33}
\end{psmatrix}
$$ (involving the aforementioned identity $\bot\tensor P\to
\bot\tensor P$) becomes:
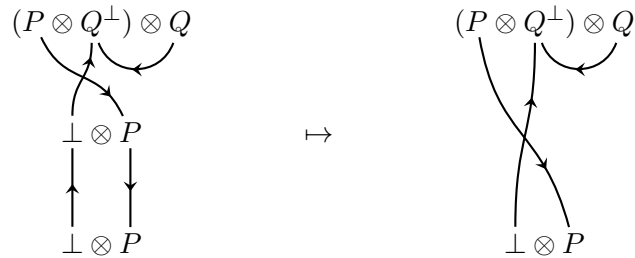
\begin{center}\begin{math}\psset{rowsep=1cm}\begin{psmatrix}
(\rnode{P}{P}\tensor\rnode{Qp}{Q}\perp)\tensor\rnode{Q}{Q} &&
(\rnode{P'}{P}\tensor\rnode{Qp'}{Q}\perp)\tensor\rnode{Q'}{Q}
\\
\nbot{b1}\tensor\rnode{P1}{P}
&
\mapsto \\
\nbot{b2}\tensor\rnode{P2}{P}
&&
\nbot{b2'}\tensor\rnode{P2'}{P}
\uvec{b2}{b1}\dvec{P1}{P2}
\vecanglespos{b1}{Qp}{90}{-90}{.8}
\vecanglespos{P}{P1}{-65}{115}{.8}
\dloopvecanglesheight{Q}{Qp}{-110}{-70}{1}
\vecanglespos{b2'}{Qp'}{90}{-90}{.7}
\vecanglespos{P'}{P2'}{-80}{105}{.7}
\dloopvecanglesheight{Q'}{Qp'}{-110}{-70}{1}
\end{psmatrix}\end{math}\end{center}
A more interesting example of composition is on
page~\pageref{firstexamplecomp} of the Introduction.

\parag{Underlying GoI category} The category $\goipfunc$, the result of applying
the geometry-of-interaction or feedback
construction $\goi$ \cite{Gir89,JSV96,Abr96}
to the category $\pfunc$ of sets and partial functions, has the
following structure.  An object is a signed set $S$, whose elements we
shall call \emph{leaves} (each signed either \emph{positive} or
\emph{negative}), and a morphism $S\to T$ is a \defn{leaf function}: a partial function from
$S\posrestr+T\negrestr$ to $S\negrestr+T\posrestr$, where
$(-)\posrestr$ (resp.\ $(-)\negrestr$) restricts to positive (resp.\
negative) leaves.  For example,
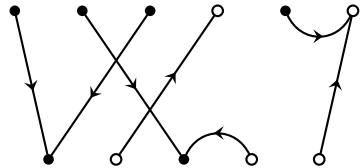
\begin{displaymath}
\psset{nodesep=-.2mm}
\begin{psmatrix}
\posb{u}\bgap\posb{uu}\bgap\posb{a}\bgap\negb{aa}\bgap\posb{uuu}\bgap\negb{uuuu} \\
\posb{a1}\bgap\negb{aa1}\bgap\posb{b1}\bgap\negb{bb1}\bgap\negb{u1}
\linevec{u}{a1}
\linevecpos{uu}{b1}{.53}
\linevecpos{a}{a1}{.59}
\linevecpos{aa1}{aa}{.59}
\uloopvecanglesheight{bb1}{b1}{115}{65}{1}
\dloopvecanglesheight{uuu}{uuuu}{-65}{-115}{1}
\linevec{u1}{uuuu}
\end{psmatrix}
\end{displaymath}
is a (total) morphism from the upper signed set (4 positive $\bullet$
and 2 negative $\circ$ leaves) to the lower one (2 positive and 3
negative leaves).
Composition is simply (finite) path composition: for an example, see
page~\pageref{goi-comp} of the Introduction.  Turbo cut elimination is
the very same path composition, hence there is a forgetful (faithful)
functor from the category $\N$ of MLL proof nets to $\goipfunc$,
extracting the leaves from a formula.  Again, see the Introduction for
examples.

\section{Linear complexity of proof net correctness}\label{linear}

\begin{theorem}[Linear complexity]
Verification of proof net correctness is linear in the number of
leaves: if $f$ is a
leaf function on a cut sequent\/ $\Gamma$, then determining whether
$f$ is a proof net can be done in linear time in the number of leaves
of\/ $\Gamma$.
\end{theorem}

\begin{proof}
Verifying the \Matching{} condition is clearly linear time.  The
\Switching{} condition is a simple corollary of the unit-free theorem
\cite{Gue99,MO00}: the function $f$ determines a standard
unit-free proof structure $\widehat{f}$ on $\widehat{\Gamma}$, as
follows.  First, replace every cut pair $\cutpair{A}$ by $A\tensor
A\perp$.  We may assume every positive leaf has an incoming $f$-edge:
every literal does, by \Matching; if the $1$ of a subformula $A\tensor
1$ doesn't, replace $A\tensor 1$ by $A$; if the $1$ of $A\parr 1$
doesn't, \Switching{} fails.
Re-label each positive literal to $1$ and each negative literal to
$\bot$.  Replace each $1$ by 
$1^n$ where $n\ge 1$ is the number of $f$-edges targetting the $1$,
and $1^n$ denotes the tensor product of $n$ copies of $1$ (bracketed
arbitrarily); re-target the $n$ edges to the old $1$ to each target a
distinct new $1$ of $1^n$.
Finally, view the symbols $\bot$ and $1$ as complementary literals, so
we have formed a standard proof structure $\widehat{f}$ on a cut-free,
unit-free MLL sequent $\widehat{\Gamma}$.  To clarify, here is
$\widehat{f}$ for $f$ the proof net on page~\pageref{example}:
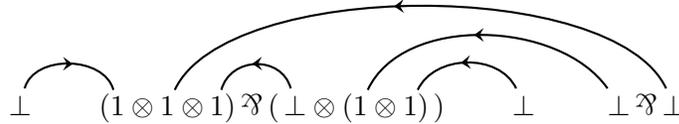
$$\begin{psmatrix}
\nbot{b1}\sgap(\one{P1}\tensor\one{P1'}\tensor\one{P1''})\,\parr\,(\,\nbot{P2}\tensor(\one{i}\tensor\one{i'})\,)\sgap\nbot{b2}\sgap\,\nbot{b3}\,\parr\,\nbot{b4}
\rule{0pt}{1.5cm}
\uloopvecheight{b1}{P1}{.8}
\uloopvecanglesheight{P2}{P1''}{105}{75}{1}
\vecanglesheight{b2}{i'}{115}{65}{.8}
\vecanglesheight{b3}{i}{125}{65}{.7}
\vecanglesheight{b4}{P1'}{115}{65}{.5}
\end{psmatrix}$$
By construction the original $f$ on $\Gamma$ is correct iff
$\widehat{f}$ on $\widehat{\Gamma}$ is correct in the usual unit-free
sense.  The construction of $\widehat{f}$ is linear time in the number
of leaves.
\end{proof}
\begin{corollary}
The theorem above extends to old proof nets (\ie, when $f$ is a
function from negative leaves to vertices of\/ $\Gamma$, optionally
with a differentiation between axiom links
\raisebox{-.3ex}{\footnotesize\,$\nbot{b}\;\;\;\one{i}\ncbar[nodesep=1pt,angleA=90,angleB=90,arm=3pt]{b}{i}$\,}
and jumps
\raisebox{-.3ex}{\footnotesize\,$\nbot{b}\;\;\;\;\;\one{i}{\psset{nodesepA=0pt,nodesepB=1pt}\uloopvecanglesheight{b}{i}{50}{130}{.9}}$\,}).
\end{corollary}
\begin{proof}
First, if differentiating, replace every axiom link
\raisebox{-.3ex}{\footnotesize\,$\nbot{b}\;\;\;\one{i}\ncbar[nodesep=1pt,angleA=90,angleB=90,arm=3pt]{b}{i}$\,}
by a jump
\raisebox{-.3ex}{\footnotesize\,$\nbot{b}\;\;\;\;\;\one{i}{\psset{nodesepA=0pt,nodesepB=1pt}\uloopvecanglesheight{b}{i}{50}{130}{.9}}$\,}.
Rewrite every compound subformula or negative leaf $A$ targeted by a
$\bot$-jump to $A\tensor 1$, and shift any $\bot$-jumps which
targeted $A$ to target the new $1$ instead.  This yields a function
$\tilde{f}$ from negative leaves to positive leaves which is correct
iff $f$ is correct; apply the above theorem to $\tilde{f}$.
To clarify, here is $\tilde{f}$ for the old proof net $f$ drawn on page~\pageref{oldnet}:
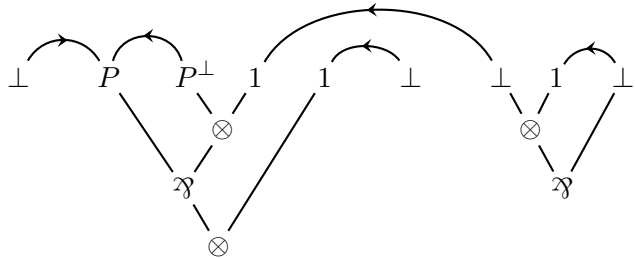
\begin{center}\begin{math}\begin{psmatrix}
\nbot{b1}\sgap\rnode{P1}{P}\,\;\;\;\;\;\,\,\rnode{P2}{P}\perp\;\;\;\;\one{i}\,\;\;\;\;\;\,\,
\one{i'}\sgap\nbot{b2}\sgap\nbot{b3}\,\;\;\;\;\one{i''}\,\;\;\;\;\:\,\nbot{b4}\rule{0pt}{1cm} \\
\parsenodeone{P2}{i}{tensor}{\tensor}
\parsenodetwoleft{P1}{i}{par}{\parr}{tensor}
\parsenodethreeright{P1}{i'}{tensor2}{\tensor}{par}
\parsenodeone{b3}{i''}{tensor3}{\tensor}
\parsenodetworight{b3}{b4}{par2}{\parr}{tensor3}
\vecanglesheight{b2}{i'}{125}{65}{.7}
\uloopvecanglesheight{P2}{P1}{105}{75}{1}
\vecanglesheight{b3}{i}{120}{60}{.7}
\uloopvecanglesheight{b1}{P1}{60}{120}{.9}
\psset{nodesep=1.5pt}
\uloopvecanglesheight{b4}{i''}{120}{60}{.9}
\end{psmatrix}\end{math}\vspace{1ex}\end{center}
The construction $f\mapsto\tilde{f}$ is
linear time in the number of leaves.
\end{proof}

\section{Laminated GoI composition for MALL nets}\label{mall}

Continuing the $\goi$ theme, we observe that composition (turbo cut
elimination) of MALL proof nets \cite{HG03,HG05} can be viewed as occurring
in a `laminated' variant of $\goipfunc$.\footnote{Since the
MALL proof nets are unit-free, we could work with the category of sets
and partial injections instead of $\pfunc$.}

Leaf functions $f:S\to T$ and $g:T\to U$ \defn{synchronise} in $T$,
denoted $f\sync_T g$ or simply $f\sync g$, if
\begin{eqnarray*}
\\[-6ex]
(\im f)\cap T &=& (\dom g)\cap T \\
(\dom f)\cap T &=& (\im g)\cap T
\end{eqnarray*}
In other words, for every (signed) leaf $t\in T$ there is an edge into
$t$ iff there is an edge out of $t$ (in $f\cup g$, viewed as a
directed graph on $S+T+U$).

Let $\finsetp$ be the category of finite sets and partial functions (a full subcategory of $\pfunc$).
Thus $\goifinsetp$ is the full subcategory of $\goipfunc$ with finite
objects (signed sets).
Finiteness will be critical for associativity of composition in the
construction below.

Define the category $\lam(\goifinsetp)$
as follows.  Objects are inherited from $\goifinsetp$ (finite signed
sets), and a morphism $S\to T$ is a set $l$ of leaf functions $S\to
T$, which we call a \defn{laminated leaf function} from $S$ to $T$.
Given laminated leaf functions $l:S\to T$ and $m:T\to U$, their
composite $\,l\synccomp m\,:\,S\to U\,$ is by
\defn{synchronised composition}:
$$l\synccomp m\;\;\;=\;\;\;\{\;f;g\text{ such that }f\in l,\;g\in
m,\text{ and }f\sync g\;\},$$ where $f;g:S\to U$ is the (sequential)
composite of leaf functions $f:S\to T$ and $g:T\to U$ in $\goifinsetp$.
Thus $l\synccomp m$ collects the pairwise composition of all synchronising
leaf functions $f\in l$ and $g\in m$.
The identity laminated leaf function $\id_S:S\to S$
comprises every leaf function $S\to S$ contained in the identity leaf
function $S\to S$ in $\goifinsetp$.
The identity law $\id_S\synccomp l=l=l\synccomp \id_S$ holds since
every leaf function $f\in l$ synchronises with a unique leaf function
in $\id_S$.
Composition is associative since given leaf functions $f:S\to T$ and
$g:T\to U$, 
$$f\sync g\;\;\;\;\;\implies\;\;\;\;(\dom f;g)\cap S\;=\;(\dom f)\cap
S\;\;\;\text{ \it and }\;\;\;(\im f;g)\cap U\;=\; (\im g)\cap U\,.$$
This \defn{stable domain-image} property depends critically on
finiteness:
if infinite signed sets are available, the property fails, and so does
associativity of composition.  Let $f:S\to T$, $g:T\to U$ and $h:U\to
V$ be the following composable leaf functions
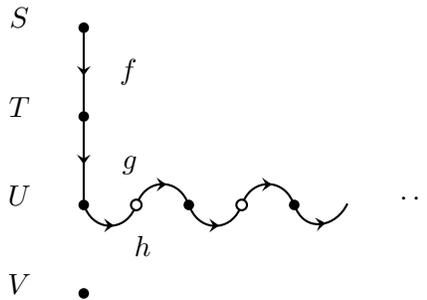
\begin{center}\vspace{1ex}\begin{math}
\newcommand{\ul}[2]{\uloopvecanglesheight{#1}{#2}{65}{115}{1}}
\newcommand{\dl}[2]{\dloopvecanglesheight{#1}{#2}{-65}{-115}{1}}
\psset{rowsep=.7cm,colsep=.7cm,nodesep=-.2mm}\begin{psmatrix}
S & \posb{s} \\
T & \posb{t} \\
U & \posb{u1} & \negb{u1'} & \posb{u2} & \negb{u2'} & \posb{u3} & \pnode{u3'} & \cdots\\
V & \posb{v}
\dvec{s}{t}\naput{\hspace{2ex}f}
\dvec{t}{u1}
\dl{u1}{u1'}\nbput{\hspace{5ex}h}
\ul{u1'}{u2}\naput{g\hspace{5ex}}
\dl{u2}{u2'}
\ul{u2'}{u3}
\dl{u3}{u3'}
\end{psmatrix}\end{math}\vspace{1ex}\end{center}
so $S,T,V$ are singletons, $U$ is countably inifinite, and $g$ and $h$
continue as suggested by the ellipsis.  Note that $f\sync g\sync h$.
The stable domain property fails for $g$ and $h$: $(\dom g)\cap T=T$
but $(\dom g;h)\cap T$ is empty (since $g;h$ is empty).  Thus
$$\,f\,\not\sync\, g;h\,\text{\hspace{3ex} but \hspace{3ex}}\,f;g\,\sync\, h\,$$ so
$\{f\}\,\synccomp\,
\big(\{g\}\synccomp\{h\}\big)\;\neq\;\big(\{f\}\synccomp\{g\}\big)\,\synccomp\,\{h\}$ for the
corresponding singleton laminated leaf functions.
(Note, however, that $f;(g;h)=(f;g);h$, so lamination is important.)

The canonical embedding $\goifinsetp\hookrightarrow\lam(\goifinsetp)$
is more interesting than one might have anticipated.  It is of course
the identity on objects (finite signed sets), but not the `identity'
on leaf functions in the naive sense of mapping a leaf function
$f:S\to T$ in $\goifinsetp$ to the singleton laminated leaf function
$\{f\}:S\to T$ in $\lam(\goifinsetp)$.  Rather, the image of $f$
comprises every leaf function $S\to T$ which is contained in $f$.
Thus the embedding takes the downard closure of $f$.  Observe how this
preserves identities.

\parag{Forgetful functor from MALL nets}  
There is a forgetful functor from the category of MALL proof nets
\cite[\S5.2]{HG05} to $\lam(\goifinsetp)$, extracting the (signed)
leaves of a MALL formula: the synchronisation property $\sync$ is the
\emph{matching} property \cite[\S5.11]{HG05}, and synchronised
composition corresponds to the definition of the \emph{normal form} of
a set of MALL linkings, by turbo cut elimination \cite[\S5.11]{HG05}
(the path composition \emph{reduction} of \cite[\S5.11]{HG05} being
composition in $\goifinsetp$).

One would hope to be able to find relationships between MALL proof
nets, $\lam(\goifinsetp)$ and the geometry of interaction for additives
as in \cite{Gir95,AJ94}.  Since MLL units are the focus of the
present paper, this is best left for another occasion.

\parag{The lamination construction}  The construction of $\lam(\goifinsetp)$ 
from $\goifinsetp$ abstracts to a general construction $\lam(\C)$ on a
category $\C$ equipped with a suitable synchronisation relation
$\sync$ between homsets: enumerate the properties of $\sync$ used
above to ensure that $\lam(\C)$ is a category.  However, whether or
not the lamination construction leads anywhere interesting remains to
be seen.  It may be useful in constructing models for the addtives,
perhaps in concert with double glueing \cite{Tan97,HS03}.

\section{Previous approaches}

Girard's passage quoted on the first page of the Introduction gives a
convenient summary of old proof nets.  Normalisation is hampered by
having to move targets of $\bot$-jumps.

Proof nets for MLL with units are given in \cite{BCST96} and
\cite{LS04}.  Neither solves the problems in Girard's quote:
each suffers from the need to move $\bot$-jumps during elimination, so
one is lumbered once again with equivalence classes.

\subsection{Circuit nets}\label{Girard-BCST}

The cut-free one-sided MLL proof nets in \cite{BCST96}
are\footnote{See the introduction to Section~2 of \cite{BCST96}.}
cut-free old proof nets (as described in Girard's quote,
page~\pageref{Gir-quote}) in circuit/wire notation, with an additional
ordering on $\bot$-jumps.  For example, the old proof net on
page~\pageref{oldnetold} is drawn thus:
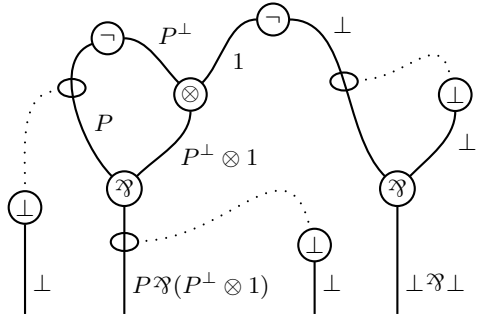
\begin{center}\footnotesize\psset{xunit=1.1cm,yunit=.5cm}
\newcommand{\botlink}[3]{\cnodeput[framesep=1pt](#1,#2){#3}{$\bot$}}
\newcommand{\conc}{1}
\newcommand{\rowtwo}{4.5}
\newcommand{\rowone}{7}
\newcommand{\paronex}{-1.3}
\newcommand{\conclabel}[2]{\nput{0}{#1}{\raisebox{7ex}{$#2$}}}
\newcommand{\wire}[4]{\nccurve[nodesep=-.2pt,angleA=#3,angleB=#4]{#1}{#2}}
\newcommand{\downwire}[2]{\ncline[nodesepA=-.2pt,angleA=-90,angleB=90]{#1}{#2}}
\begin{pspicture}(0,1.3)(0,8.7)
\cnodeput[framesep=2pt](-1.5,8.5){neg}{$\neg$}
\botlink{-2.5}{4}{b1}
\pnode(0,8){i}
\botlink{1}{3}{b2}
\cnodeput[framesep=2pt](.5,9){neg2}{$\neg$}
\pnode(1.5,8){b3}
\botlink{2.7}{\rowone}{b4}
\pnode(-2.5,\conc){c1}
\pnode(\paronex,\conc){c2}
\pnode(1,\conc){c3}
\pnode(2,\conc){c4}
\cnodeput[framesep=1pt](\paronex,\rowtwo){par1}{$\footnoteparr$}
\cnodeput[framesep=1pt](2,\rowtwo){par2}{$\footnoteparr$}
\cnodeput[framesep=1pt](-.5,\rowone){tensor}{$\tensor$}
\downwire{b1}{c1}
\downwire{par1}{c2}
\ncput[npos=.33]{\ovalnode{j2}{\;}}
\downwire{b2}{c3}
\downwire{par2}{c4}
\wire{par1}{neg}{135}{-160}
\ncput[npos=.6]{\ovalnode{j1}{\;}}
\nbput[npos=.33]{$P$}
\wire{par1}{tensor}{45}{-90}
\wire{tensor}{neg}{135}{0}
\nbput[npos=.55]{$P\perp$}
\wire{tensor}{neg2}{45}{180}
\nbput[npos=.4]{$1$}
\wire{par2}{b4}{45}{-90}
\wire{neg2}{par2}{0}{135}
\naput[npos=.25]{$\bot$}
\ncput[npos=.5]{\ovalnode{j3}{\;}}
\conclabel{c1}{\bot}
\conclabel{c2}{P\footnoteparr (P\perp\tensor 1)}
\conclabel{c3}{\bot}
\conclabel{c4}{\bot\footnoteparr\bot}
\nput{-65}{tensor}{\raisebox{-4.5ex}{$P\perp\tensor 1$}}
\nput{-70}{b4}{\raisebox{-3ex}{$\bot$}}
\psset{linestyle=dotted}
\nccurve[angleA=-170,angleB=90]{j1}{b1}
\nccurve[angleA=3,angleB=105,ncurv=1]{j2}{b2}
\nccurve[angleA=5,angleB=105,ncurv=1]{j3}{b4}
\end{pspicture}
\end{center}
Links are drawn as circular nodes, formulas are drawn as (labelled)
wires, and $\bot$-jumps are drawn dotted.  By an \emph{MLL proof net}
in the \cite{BCST96} setting we mean the special case when the base is
a set of propositional variables, and $(-)\perp$ is restricted to
propositional variables (as usual with MLL formulas).  The primary net
definition in \cite{BCST96} is two-sided; a one-sided net is simply a
two-sided net with the tensor unit $1$ on the input side (see the
paragraph following Corollary~5.3 of
\cite{BCST96}).  In drawing the one-sided net above, we omitted this
input unit and its jump.
The minor difference with old proof nets is that when multiple
$\bot$-jumps target the same wire, they are ordered along the wire; in
an old proof net there is no such ordering on $\bot$-jumps targetting
the same subformula.

The problem with normalisation (see Girard's passage on
page~\pageref{Gir-quote}) remains.  For example, if we cut against the
$P\,\parr(P\perp\tensor 1)$ wire above, we do not have a cut redex:
first we must re-wire the incoming $\bot$-jump to elsewhere in the
empire of the $\bot$; we're once again resorting to equivalence
classes for normalisation.

A key feature of the approach in \cite{BCST96} is the modularity over
negation and planarity.  Circuit nets modulo equivalence describe the
free linearly distributive and star-autonomous categories over a
polygraph (\eg, over a category), yielding full coherence.
For an internal language presentation of free star-autonomous
categories, with full coherence, see \cite{KO99} (again modulo
an equivalence/congruence).

\subsection{Syntactic nets}\label{LS}

The paper \cite{LS04} defines a proof net on a cut sequent
$\Gamma$ as a separate MLL formula $\Theta$ whose leaves from
left-to-right are a permutation of those of $\Gamma$.
The formula $\Theta$ is shown upside down above the sequent, and the
permutation is represented by permitting argument edges to cross in
the upper half.
The $\bot$-attachments and axiom links are thus enveloped in an
additional syntactic layer $\Theta$, with $\bot$-attachments as
\raisebox{-.3ex}{\psset{nodesep=1pt}\footnotesize$\rnode{A}{\;\strut}\;\;\;\;\rput(0,.3){\rnode{t}{\tensor}}\;\;\;\;
\rnode{AA}{\bot}\ncline{A}{t}\ncline{t}{AA}$}
and axiom links as
\raisebox{-.3ex}{\psset{nodesep=1pt}\footnotesize$\rnode{A}{A}\;\;\;\;\rput(0,.3){\rnode{t}{\tensor}}\;\;\;\;
\rnode{AA}{A}\perp\ncline{A}{t}\ncline{t}{AA}$}.
Here is an example of a proof net on the three-formula sequent $\bot,\,1\tensor
P,\,\bot\tensor((P\perp\tensor P\perp)\parr P)$, essentially Figure 2
of \cite{LS04}:
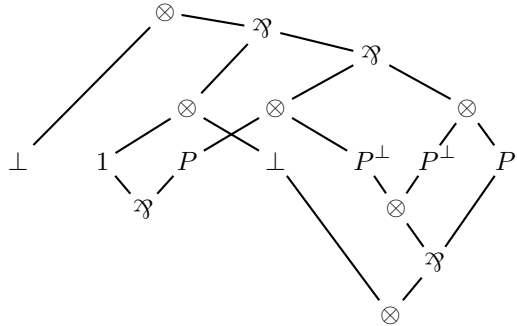
\begin{center}\small\begin{math}
\begin{psmatrix}
\nbot{b1}\sgap\one{i}\sgap\rnode{P}{P}\sgap\nbot{b2}\sgap\rnode{Pp}{P}\perp\hspace{-1ex}
\,\;\;\;\;\,\rnode{Pp2}{P}\perp\,\;\;\;\;\,\rnode{P2}{P}\,\rule{0pt}{2.1cm} \\
\parsenodeone{i}{P}{par1}{\smallparr}
\parsenodeone{Pp}{Pp2}{tensor1}{\tensor}
\parsenodetworight{Pp}{P2}{par2}{\smallparr}{tensor1}
\parsenodethreeleft{b2}{P2}{tensor2}{\tensor}{par2}
\renewcommand{\parsenodeparents}[7]{\ncline[linestyle=none]{#1}{#2}\naput[labelsep=#3]{\rnode{#4}{#5}}%
\ncline{#6}{#4}\ncline{#7}{#4}}
\renewcommand{\parseone}{3.5ex}
\renewcommand{\parsetwo}{7.3ex}
\parsenodeone{i}{b2}{tensor1'}{\tensor}
\parsenodeone{P}{Pp}{tensor2'}{\tensor}
\parsenodeone{Pp2}{P2}{tensor3'}{\tensor}
\parsenode{tensor2'}{tensor3'}{3.5ex}{par1'}{\smallparr}
\parsenodeone{tensor1'}{par1'}{par2'}{\smallparr}
\parsenodetwoleft{b1}{par1'}{tensor4'}{\tensor}{par2'}
\end{psmatrix}\end{math}\vspace*{1ex}\end{center}
As with \cite{BCST96} nets, the problem with normalisation (see
Girard's passage on page~\pageref{Gir-quote}) remains.  For
example\footnote{This example is drawn just after Lemma~4.2 in
\cite{LS04}, with different literal labels.}, if $\Gamma$ is the cut sequent
$P\perp,\cutpair{P},P\tensor Q,Q\perp,\bot$ and $\Theta$ is the proof
net given by the MLL formula $(P\perp\tensor
P)\parr\big(((P\perp\tensor P)\parr(Q\tensor Q\perp))\tensor
\bot\big)$ (with identity permutation on leaves) then the cut cannot
be reduced immediately.  First one must apply invertible linear
distributivity / commutativity / associativity to $\Theta$, subject to
the constraint of not breaking the correctness criterion (\ie, a form
of empire-rewiring \cite{Tri94,BCST96}).  Thus one is again
resorting to equivalence classes for normalisation (see Theorem~4.3 of
\cite{LS04}).
Syntactic nets modulo equivalence describe the free star-autonomous
category with strict double involution $A=A\perpp$ generated by a set.

\small
\bibliographystyle{alpha}
\bibliography{../../bib/main}
\end{document}